\documentclass[leqno]{article}%
\usepackage{amsfonts,amsmath,amsthm,amssymb,times}
\usepackage{color}
\usepackage{url}
\usepackage{graphicx,subfigure,color,float,enumerate}
\usepackage[utf8]{inputenc}
\usepackage{mathptmx}
\usepackage{graphicx}
\usepackage{caption}

\usepackage[11pt]{moresize}
\pdfoutput=1
\newtheorem{theorem}{THEOREM}[section]

\newtheorem{definition}{DEFINITION}[section]

\theoremstyle{remark}

\numberwithin{equation}{section}

\title{\bf Two Sample Test for Extrinsic Antimeans on Planar Kendall Shape Spaces with  an Application to Medical Imaging}
  \author{Aaid Algahtani    \footnote{Department of Statistics,
    Florida State University, Tallahassee, FL 32304, U.S.A.,King Saud University, Riyadh 11451, Saudi Arabia} \quad
    Vic Patrangenaru    \footnote{Department of Statistics,
    Florida State University, Tallahassee, FL 32304, U.S.A.}
  }

\begin{document}
\date{\today}
\maketitle
{\bf Abstract.}
      In this paper one develops nonparametric inference procedures for comparing two extrinsic antimeans on compact manifolds. Based on recent Central limit theorems for  extrinsic sample antimeans w.r.t. an arbitrary embedding of a compact manifold in a Euclidean space, one derives an asymptotic chi square test for the equality of two extrinsic antimeans. Applications are given to distributions on complex projective space $CP^{k-2}$ w.r.t. the Veronese-Whitney embedding, that is a submanifold representation for
      the Kendall planar shape space. Two medical imaging analysis applications are also given.

\section{Introduction}

To date, the statistical analysis on object spaces is concerned with data analysis on sample spaces that have a complete metric space structure. In include, but not limited to, data extracted from DNA and RNA sequences, from medical images (such as classification and discrimination two different groups) and 3D computer vision outputs.

In those areas of research, one needs to extract certain geometric information, from images, or from DNA.  In case of images, the features  extracted are, often times, various sorts of {\em shapes of labeled configurations of points}. For similarity shape , we refer to Kendall (1984) \cite{Ke:1984}, Kendall et al(1999) \cite{Ke:1999}, Bookstein \cite{Bo:1991}, Dryden and Mardia (2016)\cite{DrMa:2016}. For affine shapes, projective shapes, or, in general, $\mathcal G$-shapes, see Patrangenaru and Ellingson(2015), , Section 3.5 \cite{PaEl:2015}.

The development of novel statistical  principles methods for complex data types extracted from image data, relies upon key features of object spaces; namely non-linearity  and topological Structure. Information extracted from an image is typically represented as  point on non-linear spaces. Statistics and Geometry share an important commonality: they are both based on the concept of distance. However not just any distance is suitable for modeling a given object data:  distances that give a homogeneous space structure are preferred, to allow for comparing location parameters of two distribution. In addition, in general one would like for the estimators of a random object's location parameters to be consistent and furthermore allow for estimation of location parameters based on CLM results.

In classical nonparametric statistical theory on an object space $\mathcal M$ with manifold structure, an estimator's asymptotic behavior is described in terms of the asymptotic expansion around zero, of the tangential component of a consistent estimator aaround the true parameter  (see Patrangenaru and Ellingson (2015, Cha.4)\cite{PaEl:2015}), and our paper follows this direction. Fast computational algorithms are also key to the object data analysis, therefore we use an extrinsic analysis (see Bhattacharya et al.(2012)\cite{BhElLiPaCr:2012}), which, for such reasons, is preferred in Object Data Analysis, over so called ``intrinsic" approach, if there is a choice for a distance on $\mathcal M.$ The extrinsic data analysis on $\mathcal M$ is using a chord distance, which is the distance $\rho_j$ induced by the Euclidean distance via the embedding $j$ of $\mathcal M$, and is given by $\rho_j (p,q) = {\|j(q)-j(p) \|}^2_0.$ Here $\|\cdot \|_0$ is the Euclidean norm.

If for a given probability measure $Q$ on an embedded compact metric space $(\mathcal M, \rho_j),$ the Fr\'echet function associated with the random object $X,$ with $Q=P_X$ is
\begin{equation}\label{eq:frechet-function}
\mathcal F(p) =  \mathbb E( \rho_j^2 (p,X))=\int_{\mathcal{M}}  \| j(x) - j(p) \|^{2}_{0} Q(dx),
\end{equation}
and its minimizers form the extrinsic mean set of $X$ (see Bhattacharya and Patrangenaru 2003 \cite{BhPa:2003}). If there is a unique point in the extrinsic mean set, this point is the {\em extrinsic mean}, and is labeled $\mu_{j,E}(Q)$ . The maximizers in \eqref{eq:frechet-function} form the extrinsic antimean set. In there is a unique point in the extrinsic antimean set, this point is the extrisic antimean, and is labeled $\alpha\mu_{j,E}(Q)$, or simply $\alpha\mu_E,$ when $j$ is known (See Patrangenaru and Guo 2016\cite{PaYaGu:2016}).
Also, given $X_1, \dots, X_n$ i.i.d.r.o.'s from $Q$, their extrinsic sample mean $\bar X_E,$ is the extrinsic mean of the empirical distribution $\hat Q_n = {\frac{1}{n}} \sum_{i=1}^n \delta_{X_i}$ and their extrinsic sample antimean $a\bar X_E$ is the extrinsic antimean of the empirical distribution $\hat Q_n = {\frac{1}{n}} \sum_{i=1}^n \delta_{X_i},$ assuming these sample extrinsic means or antimeans exist.

The main goal of this article is to establish some general methods for comparing extrinsic antimeans on a compact manifold embedded in an Euclidean space. In Section 2, some preliminary result, a central limit theorem  for sample extrinsic antimeans on compact manifolds is brought up for future developments. Thereafter, a two sample test statistic for extrinsic antimeans on compact manifolds is given in Section 3. In Section 4, one focuses on the case of the planar Kendall shape space of {k}-ads, $\Sigma_2^k,$ which is often represented in literature by the complex projective space $\mathbb C P^{k-2}$ ( see Patrangenaru and Ellingson (2015)\cite{PaEl:2015}, Ch.2, p.142 ). The embedding considered here is the Veronese-Whitney embedding of $\mathbb C P^q$ in the space of self adjoint $(q+1) \times (q+1)$ complex matrices (see Patrangenaru and Ellingson (2015)\cite{PaEl:2015}, Ch.2,p. 154). In Section 5, we derive a nonparametric two sample test for the equality of the population Veronese-Whitey (VW) antimeans on the complex projective space $\mathbb{C}P^{k-2}$. Finally, Section 6 illustrate a two example of Kendall shape data from Bookstein ((1997)\cite{Bo:1997}). The first example is on Apert syndrome vs clinically normal children, known as the University School data. The second example is on comparing brain scan shapes in schizophrenic vs clinically normal children.

\section{Extrinsic mean and antimean of a random object}
As discussed in the introduction, when $j$ is an embedding of a compact object space into a numerical space, the set of minimizers of the Fr\'echet function is the extrinsic mean set, and when this set has a unique point, this point is called the {\em extrinsic mean}. The maximizers of this function form the extrinsic antimean set, when this later set has one point only, this point is called the extrinsic antimean.
We have following properties of the extrinsic antimean from Patrangenaru, Yao and Guo \cite{PaYaGu:2016}

\begin{definition}\label{focal}
 A point $y \in R^N $for which there is a unique point
$p\in \mathcal M$ satisfying the equality,
\begin{equation}
\sup_{x\in \mathcal M}||y-j(x)||_0= {\|y-j(p) \|}_0
\end{equation}
is called $\alpha j$-nonfocal. A point which is not $\alpha j$-nonfocal
is said to be $\alpha j$-focal. If y is an $\alpha j$-nonfocal point, its farthest projection on
$j(\mathcal M)$ is the unique point $z = P_{F,j}(y) \in  j(\mathcal M)$
with $\sup_{x\in \mathcal M}||y-j(x)||_0=\,d_0(y, j(p))$.
\end{definition}

\begin{definition}\label{focal2}
 A point $y \in R^N $for which there is a unique point
$p\in \mathcal M$ satisfying the equality,
\begin{equation}
\inf_{x \in \mathcal M}||y-j(x)||_0= {\|y-j(p) \|}_0
\end{equation}
is called $ j$-nonfocal. A point which is not $ j$-nonfocal
is said to be $ j$-focal. If y is an $ j$-nonfocal point, its projection on
$j(\mathcal M)$ is the unique point $z = P_{j}(y) \in  j(\mathcal M)$
with $\inf_{x\in \mathcal M}||y-j(x)||_0=\,d_0(y, j(p))$.
\end{definition}

\begin{definition}
A probability distribution $Q$ on $\mathcal M$
is said to be $\alpha j$-nonfocal if the mean $\mu$ of $j(Q)$ is $\alpha j$-nonfocal.

A probability distribution $Q$ on $\mathcal M$
is said to be $ j$-nonfocal if the mean $\mu$ of $j(Q)$ is $ j$-nonfocal.
\end{definition}

Then we have the following theorem from Patrangenaru, Guo and Yao (2016)\cite{PG16}, which in particular is valid for a probability measure on an embedded compact object space with the chord distance $\rho_j:$
\begin{theorem}\label{frec} {\it Let $Q=P_X$ be a probability measure  associated with the random object $X$ on a compact metric space $(M, \rho).$ So we have $F(p) = E(\rho ^2(p,X)) $ is finite on $M.$} (a)~{\it Then, given any $\varepsilon>0$, there exist a $P$-null set $N$ and $n(\omega)<\infty$ $\forall\, \omega \in N^c$ such that the Fr\'echet (sample) antimean set of $\hat Q_n =\hat Q_{n,\omega}$ is contained in the $\varepsilon$-neighborhood of the Fr\'echet antimean set of $Q$ for all $n\geq n(\omega)$.} (b)~{\it If the Fr\'echet antimean of $Q$ exists then every measurable choice from the Fr\'echet (sample) antimean set of $\hat Q_n$ is a strongly consistent estimator of the Fr\'echet antimean of} $Q$.
\end{theorem}

\subsection{Previous Asymptotic Results for Extrinsic Sample Antimeans }
In preparation, we are using the large sample distribution for extrinsic sample antimeans given in Patrangenaru et al (2016 \cite{PaYaGu:2016}).

Assume $j$ is an embedding of a $d$-dimensional manifold $\mathcal M$ such that $j(\mathcal M)$ is closed in $\mathbb R^N$, and $Q = P_X$ is a $\alpha j$-nonfocal probability measure on $\mathcal M$ such that $j(Q)$ has finite moments of order 2. Let $\mu$ and $\Sigma$ be the mean and covariance matrix of $j(Q)$ regarded as a probability measure on $\mathbb R ^N$. Let $\mathcal F$ be the set of $\alpha j$-focal points of $j(\mathcal M)$, and let $P_{F,j}:{\mathcal F}^c \to j(\mathcal M)$ be the farthest projection on $j(\mathcal M)$. $P_{F,j}$ is differentiable at $\mu$ and has the differentiability class of $j(\mathcal M)$ around any $\alpha j$ nonfocal point. In order to evaluate the differential $d_{\mu}P_{F,j}$ we consider a special orthonormal frame field that will ease the computations.

A local frame field $p \to (e_1(p), \dots, e_k(p))$, defined on an open neighborhood $U \subseteq \mathbb R^N$ is {\em adapted to the embedding $j$} if it is an orhonormal frame field and $\forall x \in j^{-1}(U), e_r(j(x)) = d_x j(f_r(x)), r\in\{1,\ldots, d\},$ where $(f_1, \dots, f_d)$ is a local frame field on $\mathcal M,$ and $f_r(x)$ is the value of the local vector field $f_r$ at $x.$

Let $e_1, \dots, e_N$ be the canonical basis of $\mathbb R^N$ and assume $(e_1(p), \dots, e_N(p))$ is an adapted frame field around $P_{F,j} (\mu) =  j(\alpha \mu_E)$. Then $d_\mu P_{F,j} (e_b) \in T_{P_{F,j}(\mu)}j(\mathcal M)$ is a linear combination of $e_1(P_{F,j} (\mu)), \dots, e_d(P_{F,j} (\mu))$:
\begin{equation}
    d_\mu P_{F,j} (e_b) = \sum^d_{a=1} (d_\mu P_{F,j} (e_b)) \cdot e_a(P_{F,j} (\mu))e_a(P_{F,j} (\mu))
\end{equation}
{where $d_\mu P_{F,j}$ is the differential of $P_{F,j}$ at $\mu.$} By the delta method, $n^{1/2}(P_{F,j} (\overline{j(X)}) - P_{F,j}(\mu))$ converges weakly to $N_N(0_N,\alpha \Sigma_\mu)$, where $\overline{j(X)} = \frac{1}{n} \sum^n_{i=1} j(X_i)$ and
\begin{equation}
    \begin{aligned}
       \alpha \Sigma_\mu = [\sum^d_{a=1} d_\mu P_{F,j} (e_b) \cdot e_a(P_{F,j} (\mu))e_a(P_{F,j} (\mu))]_{b\in\{1,\ldots, N\}} \\ \times \Sigma [\sum^d_{a=1} d_\mu P_{F,j} (e_b) \cdot e_a(P_{F,j} (\mu))e_a(P_{F,j} (\mu))]^T_{b\in\{1,\dots,N\}}{.}
    \end{aligned}
\end{equation}
{ Here} $\Sigma$ is the covariance matrix of $j(X_1)$ w.r.t the canonical basis $e_1, \dots, e_N.$

The asymptotic distribution $N_N(0_N,\alpha \Sigma_\mu)$ is degenerate and the support of this distribution is on $T_{P_{F,j}}j(\mathcal M)$, since the range of $d_\mu P_{F,j}$ is $T_{P_{F,j(\mu)}}j(\mathcal M)$. Note that $d_\mu P_{F,j} (e_b) \cdot e_a(P_{F,j} (\mu)) = 0$ for $a\in\{d+1, \dots,N\}$.

The tangential component $tan(v)$ of $v \in \mathbb R^N$, w.r.t the basis $e_a(P_{F,j}(\mu)) \in T_{P_{F,j(\mu)}}j(\mathcal M), a \in\{1, \dots, d\}$ is given by
\begin{equation}
    tan(v) = [e_1{(P_{F,j} (\mu))}^T v, \dots, e_d{(P_{F,j} (\mu))}^T v]^T{.}
\end{equation}
Then, the random vector ${(d_{\alpha \mu_{E}}j)}^{-1}(tan(P_{F,j} (\overline{(j(X))})-P_{F,j}(\mu))) = \sum^d_{a=1} {\overline{X}}^a_j f_a$ has the following  anticovariance matrix w.r.t the basis $f_1(\alpha \mu _E), \dots, f_d(\alpha \mu_E)$:
\begin{equation}
    \begin{aligned}
        \alpha \Sigma_{j,E} =   e_a{(P_{F,j} (\mu))}^T \alpha \Sigma_\mu e_b{(P_{F,j} (\mu))}_{1 \leq a, b \leq d} \\ = [d_\mu P_{F,j} (e_b) \cdot e_a(P_{F,j} (\mu))]_{a \in\{1, \dots, d\}} ~ \Sigma  ~  [d_\mu P_{F,j} (e_b) \cdot e_a(P_{F,j} (\mu))]^T _{a \in\{1, \dots, d\}}{,}
    \end{aligned}
\end{equation}
which is the {\em anticovariance} matrix of the random object $X$.
To simplify the notation, we set
\begin{equation}\label{e:b}
B=[d_\mu P_{F,j} (e_b) \cdot e_a(P_{F,j} (\mu))]_{a \in\{1, \dots, d\}}
\end{equation}
Similarly, given i.i.d.r.o.'s $X_1, \cdots, X_n$ from $Q$, we define the {\em sample anticovariance} matrix $aS_{j,E,n}$ as the anticovariance matrix associated with the empirical distribution $\hat{Q}_n.$

If, in addition, rank
$\alpha \Sigma_\mu = d$, then $\alpha \Sigma_{j,E} $ is invertible and if we define the
$j$-{\em standardized antimean vector}
\begin{equation}
\overline{Z}_{j,n}=:n^{1\over 2}{\alpha\Sigma_{j,E}}^{-\frac{1}{2}}
({\bar{X}^{1}}_j, \cdots, {\bar{X}^{d}}_j)^T,
\end{equation}
using basic large sample theory results, including a generalized Slutsky's lemma ({ see Patrangenaru and Ellingson(2015)\cite{PaEl:2015}, p.65}), one has:
\begin{theorem}\label{t:anti-clt} Assume $\{X_r\}_{r=1,...,n}$ is a
random sample from the $\alpha j$-nonfocal distribution $Q$, and let
$\mu =E(j(X_1))$. Let
$(e_1(p),e_2(p) ,....,e_N(p))$ be an orthonormal
frame field adapted to  $j$. Then (a) the tangential component at the extrinsic
antimean $\alpha\mu_E$ of ${d_{\alpha \mu_{E}}j^{-1}} tan_{P_{F,j}(\mu)}(P_{F,j}(\overline{(j(X))})-P_{F,j}(\mu))$ has asymptotically a multivariate normal
distribution in the tangent space to $M$ at $\alpha\mu_E$
with mean  $0_d$ and  anticovariance matrix $n^{-1}\alpha\Sigma_{j,E}$ ,
 (b) if $\alpha\Sigma_{j,E}$
is nonsingular, the j-standardized antimean vector $\overline{Z}_{j,n} = \alpha\Sigma_{j,E}^{-\frac{1}{2}}tan_{P_{F,j}(\mu)}(P_{F,j}(\overline{(j(X))})-P_{F,j}(\mu))$
 converges weakly to a random vector with a $N_d( 0_d,I_d)$ distribution, and (c) under the assumptions of (b)
\begin{equation}\label{eq:clt-anti}\|(aS_{j,E,n})^{-\frac{1}{2}}tan_{P_{F,j}(\mu)}(P_{F,j}(\overline{(j(X))})-P_{F,j}(\mu))\|^2 \to_d \chi^2_d.\end{equation}
\end{theorem}
\section{A two sample test for extrinsic antimeans}
We now turn to two-sample tests for extrinsic antimeans of distributions on an arbitrary $m$ dimensional compact manifold $M.$ For a large sample test for equality of two extrinsic means see Bhattacharya and Bhattacharya (2012) \cite{bhbh08}, p.42. Let {$X_{ak_a} : k_a = 1, \dots, n_a, a = 1, 2$} be two independent random samples drawn from distributions { $Q_a, a = 1, 2$ } on $M,$ and let $j$ be an embedding of $M$ into $\mathbb R^N.$  Denote by $\mu_a$ the mean of the induced probability $Q_a\circ j^{-1}$ and $\alpha\Sigma_{a,j}$  its  anticovariance matrix $(a = 1, 2).$ Then the extrinsic antimean of $Q_a$ is { $\alpha\mu_{a,j} = j^{-1}(P_{F,j}(\mu_a))$}, assuming { $Q_a$} is $\alpha$-nonfocal. Write { $Y_{ak_a} = j(X_{ak_a}) k_a = 1, \dots, n_a, a = 1, 2$ } and let $a\bar Y_a, a =1, 2$  be the corresponding sample antimeans. By Theorem \ref{t:anti-clt}
\begin{equation}\label{e:clt-extr-separat-each-samp}
\sqrt{n_a}B_a[P_{F,j}(\bar Y_a ) - P_{F,j}(\mu_a)] \to_d \mathcal N_m(0, \alpha\Sigma_{a,j}), a =1, 2,
\end{equation}
{ where $\alpha\Sigma_{a,j}$ is the extrinsic anticovariance matrix of $Q_a,$ and $B_a$ are the same as in Theorem \ref{t:anti-clt}, and \eqref{e:b}, with $Q$ replaced by $Q_a$ (a=1,2). }
That is, $B_a(y)$ is the $m\times N$ matrix of an orthonormal basis (frame) of $T_y(j(M)) \subset T_y(\mathbb R^N) = \mathbb R^N$ for $y$ in a neighborhood of $P_{F,j}(\mu_a),$ and $B_a = B_a(P_{F,j}(\mu_a)).$  Similarly,
 { $C_a \doteq \Sigma_{\mu_a, a} = (D P_{F,j})_{\mu_a} \Sigma_a (D P_{F,j})_{\mu_a} ^T,  (a = 1, 2).$  The null hypothesis $H_0: \alpha\mu_{1,j} = \alpha\mu_{2,j},$ say, is equivalent to $H_0: P_{F,j}(\mu_1) = P_{F,j}(\mu_2) = \pi,$ say. Then, under the null hypothesis, letting $B = B(\pi),$ one has $B_1 = B_2 = B,$ and
 \begin{eqnarray}\label{e:clt-extr-2-sample-difference}
 [B (\frac{1}{n_1}C_{1} + \frac{1}{n_2}C_{2})B^T]^{-1/2} B[P_{F,j}( \bar Y_1) - P_{F,j}( \bar Y_2)]\to_d \mathcal N_m(0_m, I_m), \nonumber \\ \text{as} \ n = n_1+n_2 \to \infty, \text{and} \frac{n_1}{n}\to \lambda \in (0,1).
 \end{eqnarray}}
For statistical inference one estimates $C_a$ by
\begin{equation}\label{eq:vova-2sample}
\hat C_a =  (D P_{F,j})_{\bar Y_a}\hat \Sigma_{a}(D P_{F,j})_{\bar Y_a}^T
\end{equation}
where  $\hat \Sigma_a$  is the sample covariance matrix of sample $a  ( a = 1, 2).$  Also $B$ is replaced by
 $ \hat B = B(\hat \pi )$ where $\hat \pi$  is a sample estimate of $\pi.$ Under $H_0,$ both $P_{F,j}( \bar Y_1)$ and $P_{F,j}( \bar Y_2)$ are consistent estimates of $\pi,$ so we take a ``pooled estimate"
 \begin{equation}\label{e:pooled-sampel-mean}
 \hat \pi = P_{F,j}( \frac{1}{n_1 + n_2}(n_1 P_{F,j}( \bar Y_1) + n_2 P_{F,j}( \bar Y_2) )).
    \end{equation}
 We, therefore, have the following result:

\begin{theorem}\label{t:pooled-asym-ext-two-sample}
Assume the extrinsic sample  anticovariance matrix $\hat \alpha\Sigma_{a,j}$  is nonsingular for $ a = 1, 2.$ Then, under $H_0: \alpha\mu_{1, j} = \alpha\mu_{2, j},$ one has:
\begin{eqnarray} \label{e:pooled-asym-ext-two-sample}
( \hat B[P_{F,j}(\bar Y_1) -P_{F,j}(\bar Y_2)])^T[\hat B (\frac{1}{n_1}\hat C_1 + \frac{1}{n_2}\hat C_{2})\hat B^T]^{-1}
( \hat B[P_{F,j}(\bar Y_1) -P_{F,j}(\bar Y_2)])\nonumber \\
\to_d \chi_m^2 ,\\
\nonumber \\ \text{as} \ n = n_1+n_2 \to \infty, \text{and} \frac{n_1}{n}\to \lambda \in (0,1).\nonumber
\end{eqnarray}
\end{theorem}

\section{VW Means and VW Antimeans} \label{s:5}

\subsection{VW Means and VW Antimeans on $\mathbb C P^m$}
We consider the case of a probability measure $Q$ on the complex projective space $\mathcal M = \mathbb C P^q$. If we consider the action of the multiplicative group $\mathbb C ^* = \mathbb C \backslash \{0\}$ on $\mathbb C^{q+1} \backslash \{ 0 \}$, given by scalar multiplication
\begin{equation}
    \alpha (\lambda, z) = \lambda z
\end{equation}
the quotient space is the $q$-dimensional \textit{complex projective space} $\mathbb C P^q$, set of all complex lines in $\mathbb C^{q+1}$ going through $0 \in \mathbb C^{q+1}$. One can show that $\mathbb C P^q$ is a 2$q$ dimensional real analytic manifold, using transition maps, similar to those in the case of $\mathbb R P^q$ \ (see Patrangenaru and Ellingson (2015) \cite{PE15}, chapter 3).

Here we are concerned with a landmark based non-parametric analysis of similarity shape data (see Bhattacharya and Patrangenaru (2014)\cite{BhPa:2014}). Our analysis is for extrinsic antimeans. For landmark based shape data, one considers a $k$-ad $x = (x^{1}, \dots, x^{k})  \in (\mathbb{R}^{m})^{k}$, which consists of $k$ labeled points in $\mathbb{R}^{m}$ that represent coordinates of k-labeled landmarks.

In this subsection $\mathcal{G}$ is the group direct similarities of $\mathbb{R}^{m}$. A \textit{similarity} is a function $f: \mathbb{R}^{m} \to \mathbb{R}^{m}$, that uniformly scale the Euclidean distances, that is, for which there is $K>0$, such that $\lVert f(x)-f(y)\rVert = K\lVert x-y \rVert, \forall x,y \in \mathbb{R}^{m}$. Using the fundamental theorem of Euclidean geometry, one can show that a similarity is given by $f(x) = Ax+b, A^{T}A = cI_{m},c>0$. A \textit{direct similarity} is a similarity of this form, where $A$ has a positive determinant. Under composition direct similarities from a group. The object space considered here consist in orbits of the group action of $\mathcal G$ on $k$-ads, and is called direct similarities \textit{shape space of k-ads in $\mathbb{R}^{m}$}.

For $m=2$, or $m=3$, a \textit{direct similarity (Kendall) shape} is that the geometrical information that remains when location, scale and rotational effects are filtered out from a $k$-ads. Two $k$-ads $(z_{1}, z_{2}, \dots, z_{k})$ and $(z^{'}_{1}, z^{'}_{2}, \dots, z^{'}_{k})$ are said to have the \textit{same shape} if there is a direct similarity $T$ in the plane, that is, a composition of a rotation, a translation and a homothety such that $T(z_{j}) = z^{'}_{j}$ for $j = 1, \dots, k$. Having the same Kendall shape is an equivalence relationship in the space of planar $k$-ads, and the set of all equivalence classes of nontrivial $k$-ads is called the \textit{Kendall planar shape space of $k$-ads}, which is denoted  $\Sigma_{2}^{k}$ (see Balan and Patrangenaru(2005) \cite{balan2005geometry}).

\subsection{VW Antimeans on $\mathbb C P^m$}

Kendall (1984) \cite{Ke:1984} introduced that planar direct similarity shapes of k-ads; which is a set of k labeled points at least two of which are distinct, can be represented as points on a complex projective space $\mathbb{C}P^{k-2}$.
A standard shape analysis method was also introduced by Kent(1992)\cite{Kent:1992} and is called Veronese-Whitney(VW) embedding of $\mathbb{C}P^{k-2}$ in the space of $(k-1) \times (k-1)$ self adjoint complex matrices, to represent shape data in an Euclidean space. This Veronese--Whitney embedding
$j: \mathbb{C}P^{k-2} \to S(k-1,\mathbb{C})$, where $S(k-1,\mathbb{C})$ is the space of  $(k-1) \times (k-1)$ Hermitian matrices, is given by
\begin{equation} \label{veronese}
j([{\bf z}]) = {\bf z z^*, z^*z}=1.
\end{equation}

This embedding is a $SU(k - 1)$ equivariant embedding and $SU(k - 1)$ is called the special unitary group $(k - 1) \times (k - 1)$ matrices of determinant 1, since $j([Az]) = Aj([z])A^*$,  $\forall A \in SU(k - 1)$. The corresponding extrinsic mean (set) of a random shape X on $\mathbb{C}P^{k-2}$ is called the VW mean (set) (see Patrangenaru and Ellingson (2015)\cite{PaEl:2015}, Ch.3), and the VW mean,when it exists,and is labeled $\mu_{VW}(X), \mu_{VW}$ or simply $\mu_E$.The corresponding extrinsic antimean (set) of a random shape X,is called the VW antimean (set) and is labeled $\alpha\mu_{VW}(X),\alpha\mu_{VW}$ or $\alpha\mu_E$.

For the VW-embedding of the complex projective space, we have the following theorem for sample VW-means from Bhattacharya and Patrangenaru(2003) \cite{BhPa:2003}:

\begin{theorem}
Let $Q$ be a probability distribution on $\mathbb{C}P^{k-2}$ and let $\{[Z_{r}], \parallel Z_{r} \parallel = 1, {r = 1, \dots, n}\}$ be a i.i.d.r.o.'s from $Q$. $(a)$ $Q$ is nonfocal iff $\lambda$ the largest eigenvalue of $E[Z_{1}Z_{1}^{*}]$ is simple and in this case $\mu _{E}{Q} = [m],$ where $m$ is an eigenvector of $E[Z_{1}Z_{1}^{*}]$ corresponding to $\lambda$, with $\parallel m \parallel = 1$. $(b)$ The extrinsic sample mean $\overline{X}_{E} = \overline{[m]}$, where $m$ is an eigenvector of norm 1 of $J = \frac{1}{n} \sum^n_{i=1}Z_i Z^*_i$, $\| Z_i\| = 1, i = 1, \dots, n$, corresponding to the largest eigenvalue of J.
\end{theorem}

We also have a similar result for sample VW-antimeans (see Wang, Patrangenaru and Guo \cite{WaPaGu:20}):
\begin{theorem}\label{th:c-antimean}
Let $Q$ be a probability distribution on $\mathbb{C}P^{k-2}$ and let $\{[Z_{r}], \parallel Z_{r} \parallel = 1_{r = 1, \dots, n}\}$ be a i.i.d.r.o.'s from $Q$. $(a)$ $Q$ is $\alpha $-nonfocal iff $\lambda$, the smallest eigenvalue of $E[Z_{1}Z_{1}^{*}]$ is simple and in this case $\alpha \mu _{j,E}{Q} = [m],$ where $m$ is an eigenvector of $E[Z_{1}Z_{1}^{*}]$ corresponding to $\lambda$, with $\parallel m \parallel = 1$. $(b)$ The extrinsic sample antimean $\alpha \overline{X}_{E} = {[m]}$, where $m$ is an eigenvector of norm 1 of $J = \frac{1}{n} \sum^n_{i=1}Z_i Z^*_i$, $\| Z_i\| = 1, i = 1, \dots, n$, corresponding to the smallest eigenvalue of J.
\end{theorem}

\section{Two sample testing problem for VW antimeans on $\Sigma_2^k$}

\subsection{Application on the planar shape space of $k$-ads}

We are concerned with a landmark based nonparametric analysis of similarity shape data. For landmark based shape data , we will denote a $k$-ad by $k$ complex numbers: $z_j = x_j+iy_j, ~ 1 \leq j \leq k.$ We center the $k$-ad at $\langle z \rangle = {1\over k} \sum\limits_{j=1}^{k} z_j$ via a translation; next we rotate it by an angle $\theta$ and scale it, operations that are achieved by multiplying $z- \langle z \rangle$ by the complex number $\lambda = r e^{i\theta}.$ We can represent the shape of the $k$-ad as the complex one dimensional vector subspace passing through $0$ of the linear subspace $L^{k-1} = \{\zeta \in \mathbb C^k, 1_k^T\zeta = 0\}.$  Thus, the space of $k$-ads is the set of all complex lines on the hyperplane, $L^{k-1} = {w \in {\mathbb C}^k \backslash \{0\}}:~ \sum\limits_1^k w_j =0.$ Therefore the shape space $\Sigma_2^k$ of nontrivial planar $k$-ads can be represented as the complex projective space $\mathbb CP^{k-2},$ the space of all complex lines through the origin in $\mathbb C^{k-1}$ using an isomorphism of $L^{k-1}$ and $\mathbb C^{k-1}.$ As in the case of $\mathbb CP^{k-2}$, it is convenient to represent the element of $\Sigma_2^k$ corresponding to a $k$-ad $z$ by the curve $\gamma(z) = [z] = {e^{i\theta}((z- \langle z \rangle)/\|z- \langle z \rangle \|):~ 0 \leq \theta \leq 2 \pi}$ on the unit sphere in $L^{k-1} \approx \mathbb C ^{k-1}.$

\subsection{Test for VW antimeans on planar shape spaces}
Let $Q_1$ and $Q_2$ be two probability measures on the shape space $\Sigma_2^k,$ and let $\alpha \mu_1$ and $\alpha \mu_2$ denote the antimeans of $Q_1\circ j^{-1}$ and $Q_2 \circ j^{-1},$ where $j$ is the VW-embedding that $j([z]) = u u ^*,$  where $u = {1\over{\|z\|}}z,$ thus $u^*u = 1$. Suppose $[W_1], \cdots, [W_n]$ and $[Z_1], \cdots, [Z_m]$ are i.i.d. random objects from $Q_1$ and $Q_2$ respectively. Let $X_i = j([W_i]),$ $Y_j = j ([Z_j])$ be their images in $j(\mathbb CP^{k-2})$ which are random samples from $Q_1 \circ j^{-1}$ and $Q_2 \circ j^{-1},$ respectively. Suppose we are to test if the VW antimeans of $Q_1$ and $Q_2$ are equal, i.e.
$$ H_0 : \alpha \mu_1 = \alpha \mu_2$$

It is known that $\alpha \mu_1=j^{-1}(P_{F}(v_1))$, where $v_1=E(X_1)$ and similarly, $\alpha \mu_2=j^{-1}(P_{F}(v_2))$, where $v_2=E(Y_1)$. We assume that both $ v_1$ and $ v_2$ have simple smallest eigenvalues.  Then under $H_0,$ their unit corresponding eigenvectors differ by a rotation.

Choose $v \in S(k-1, \mathbb C),$ with the same farthest projection as $v_1$ and $ v_2$. Suppose $ v = u\Lambda u^*,$ and $\Lambda = Diag(\lambda_1 < \lambda_2 \leq \cdots \leq \lambda_{k-1}),$ where $\lambda_a, a =1, \dots, k-1,$ are the eigenvalues of $v,$ and $u = [u_1, u_2, \cdots, u_{k-1}]$ are the corresponding eigenvectors. Also, we obtain an orthonormal basis for $S(k-1, \mathbb C)$, which is given by $\{\upsilon_b^a: 1 \leq a \leq b \leq k-1\}$ and $\{\omega_b^a: 1 \leq a \leq b \leq k-1\}$. Where $\upsilon_b^a$ has all entries zero except for those in the positions (a, b) and (b, a) that are equal to 1 and $\omega_b^a$ is a matrix with all entries zero except for those in the positions (a, b) and (b, a) that are equal to i, respectively -i.

It is defined as,

\begin{equation}
    \upsilon_b^a = \left\{
             \begin{array}{lr}
             \frac{1}{\sqrt 2} (e_a e_b^t+e_b e_a^t), a < b \\
             e_a e_a^t, a=b.
             \end{array}
\right.
\end{equation}

\begin{equation}
\omega_b^a = + \frac{i}{\sqrt{2}} (e_a e_b^t-e_b e_a^t), a < b
\end{equation}

where $\{e_a: 1 \leq a \leq k-1\}$  is the standard canonical basis for $\mathbb{C}^{N}$. For any $u \in SU(k-1) (uu^*=u^*u=1, det(u)=+1)$ where SU is special unitary group of all $(k-1) \times (k-1)$ complex matrices. $\{u\upsilon_b^au^*: 1 \leq a \leq b \leq k-1\}$ and $\{u\omega_b^au^*: 1 \leq a \leq b \leq k-1\}$ are an orthonormal frame for $S(k-1, \mathbb C)$. Now, we will choose the orthonormal basis frame $\{u\upsilon_b^au^*,u\omega_b^au^* \}$
for $S(k-1, \mathbb C)$. Then it can be shown that

\begin{equation}
    d_v P_F(u\upsilon_b^au^*) = \left\{
             \begin{array}{lr}
             0,  & \text{if}~ 1 \leq a \leq b < k-1, a=b=k-1, \\
             (\lambda_{k-1}-\lambda_a)^{-1} u\upsilon_b^au^*,  & \text{if}~ 1 \leq a  < k, b=k-1.
             \end{array}
\right.
\end{equation}

\begin{equation}
    d_v P_F(u\omega_b^au^*) = \left\{
             \begin{array}{lr}
             0,  & \text{if}~ 1 \leq a \leq b < k-1, a=b=k-1, \\
             (\lambda_{k-1}-\lambda_a)^{-1} u\omega_b^au^*,  & \text{if}~ 1 \leq a  < k , b=k-1.
             \end{array}
\right.
\end{equation}

Then, we write

\begin{equation}
\bar X-v= \mathop{\sum\sum}_{1 \leq a \leq b < k-1} \langle\,(\bar X-v),u\upsilon_b^au^* \rangle u\upsilon_b^au^* +\mathop{\sum\sum}_{1 \leq a \leq b < k-1} \langle\,(\bar X-v),u\omega_b^au^* \rangle u\omega_b^au^*
\end{equation}

Then from equations (5.3) ,(5.4) and (5.5), it follows that
\begin{align*}
    d_{v}P_F(\bar{X}-v) &= \sum\limits_{a=2}^{k-1}\sqrt{2}Re(u_a^*\bar{X}u_1){(\lambda_a-\lambda_1)}^{-1}u\upsilon_a^1 u^* \\
    &\quad + \sum\limits_{a=2}^{k-1}\sqrt{2}Im(u_a^*\bar{X}u_1){(\lambda_a-\lambda_1)}^{-1}u\omega_a^1 u^* \\
    &= \sum\limits_{a=2}^{k-1}{(\lambda_a - \lambda_1)}^{-1}(u_a^*\bar{X}u_1)u_au_1^* \\ &\quad + \sum\limits_{a=2}^{k-1}{(\lambda_a - \lambda_1)}^{-1}(u_1^*\bar{X}u_a)u_1u_a^*
\end{align*}

\begin{align*}
    d_{v}P_F(\bar{Y}-v) &= \sum\limits_{a=2}^{k-1}\sqrt{2}Re(u_a^*\bar{Y}u_1){(\lambda_a-\lambda_1)}^{-1}u\upsilon_a^1 u^* \\
    &\quad + \sum\limits_{a=2}^{k-1}\sqrt{2}Im(u_a^*\bar{Y}u_1){(\lambda_a-\lambda_1)}^{-1}u\omega_a^1 u^* \\
    &= \sum\limits_{a=2}^{k-1}{(\lambda_a - \lambda_1)}^{-1}(u_a^*\bar{Y}u_1)u_au_1^* \\ &\quad + \sum\limits_{a=2}^{k-1}{(\lambda_a - \lambda_1)}^{-1}(u_1^*\bar{Y}u_a)u_1u_a^*
\end{align*}

Define
\begin{equation}
    {T(\alpha \mu)}_{ij} = \left\{
             \begin{array}{lr}
             Re(u_{i+1}^*X_ju_1), & \text{if}~ 1 \leq i \leq k-2, ~1 \leq j \leq n,\\
             Im(u_{i-k+3}^* X_j u_1), & \text{if}~ k-1 \leq i \leq 2k-4,~ 1 \leq j \leq n.
             \end{array}
\right.
\end{equation}

\begin{equation}
    {S(\alpha \mu)}_{ij} = \left\{
             \begin{array}{lr}
             Re(u_{i+1}^*Y_ju_1), & \text{if}~ 1 \leq i \leq k-2, ~1 \leq j \leq m,\\
             Im(u_{i-k+3}^* Y_j u_1), & \text{if}~ k-1 \leq i \leq 2k-4,~ 1 \leq j \leq m.
             \end{array}
\right.
\end{equation}

Then we have,
\begin{equation}
    \bar{T}(\alpha \mu) = \frac{1}{n}\sum\limits_{j=1}^n T(\alpha\mu), ~ \bar{S}(\alpha \mu) = \frac{1}{m}\sum\limits_{j=1}^m S(\alpha\mu)
\end{equation}

Under $H_0$, $\bar{T}(\alpha \mu)$ and $\bar{S}(\alpha \mu)$ have mean zero, and as $n,m \to \infty$, suppose $(n/(m+n)) \to p,$ $(m/(m+n)) \to q,$ for some $p,q > 0;$ $p+q = 1.$ It follows that

\begin{equation}
    \sqrt{n}\bar{T}(\alpha \mu) \overset{\mathcal L}{\longrightarrow} N(0, \Sigma_1(\alpha \mu)), ~ \sqrt{m}\bar{S}(\alpha \mu) \overset{\mathcal L}{\longrightarrow} N(0, \Sigma_2(\alpha \mu))
\end{equation}

where $\Sigma_1(\alpha \mu)$ and $\Sigma_2(\alpha \mu)$ are the  covariances of ${T(\alpha \mu)}_{.1}$ and ${S(\alpha \mu)}_{.1},$ respectively. Then

\begin{equation}
    \sqrt{(n+m)}(\bar{T}(\alpha \mu)-\bar{S}(\alpha \mu))  \overset{\mathcal L}{\longrightarrow} N_{2k-4}(0, \frac{1}{p}\Sigma_1(\alpha \mu) + \frac{1}{q}\Sigma_2(\alpha \mu))
\end{equation}

Thus assuming $\Sigma_1(\alpha \mu),$ $\Sigma_2(\alpha \mu)$ and $\frac{1}{p}\Sigma_1(\alpha \mu) + \frac{1}{q}\Sigma_2(\alpha \mu)$ to be nonsingular,

\begin{equation}\label{e:asy-chi-diff}
    (n+m){(\bar{T}(\alpha \mu)-\bar{S}(\alpha \mu))}'{(\frac{1}{p}\Sigma_1(\alpha \mu) + \frac{1}{q}\Sigma_2(\alpha \mu))}^{-1}(\bar{T}(\alpha \mu)-\bar{S}(\alpha \mu)) \overset{\mathcal L}{\longrightarrow} \chi_{2k-4}^2
\end{equation}

We may take $j(\alpha \mu) = pj(\alpha \mu_1) + qj(\alpha \mu_2).$ Since $\alpha \mu_1$ and $\alpha \mu_2$ are unknown, we may estimate $\alpha \mu$ by the pooled sample mean $j(\hat{\alpha \mu}) = (n\bar{X} + m\bar{Y})/(m+n)$. Note that $\hat{\Sigma}_1(\hat{\alpha \mu})$ and $\hat{\Sigma}_2(\hat{\alpha \mu})$ are consistent estimator of $\Sigma_1(\alpha \mu)$ and $\Sigma_2(\alpha \mu)$. Thus, we may use $\hat{\Sigma}_1(\hat{\alpha \mu})$ and $\hat{\Sigma}_2(\hat{\alpha \mu})$ to obtain a two sample test statistic, where

\begin{equation}
    \hat{\Sigma}_1(\hat{\alpha \mu}) = \frac{1}{n}T(\alpha \mu){T(\alpha \mu)}^{'} - \bar{T}(\alpha \mu){\bar{T}(\alpha \mu)}^{'}
\end{equation}

\begin{equation}
    \hat{\Sigma}_2(\hat{\alpha \mu}) = \frac{1}{m}S(\alpha \mu){S(\alpha \mu)}^{'} - \bar{S}(\alpha \mu){\bar{S}(\alpha \mu)}^{'}
\end{equation}
Then the two sample test statistic can be estimated by
\begin{equation}\label{e:aTnm}
    aT_{nm} = {(\bar{T}(\hat{\alpha \mu})-\bar{S}(\hat{\alpha \mu}))}^{'}{(\frac{1}{n}\hat{\Sigma}_1(\hat{\alpha \mu})+\frac{1}{m}\hat{\Sigma}_2(\hat{\alpha \mu}))}^{-1}((\bar{T}(\hat{\alpha \mu})-\bar{S}(\hat{\alpha \mu})))
\end{equation}

Given the significance level $\beta$, we reject $H_0$ if
\begin{equation}\label{e:aTnm-level-beta}
    aT_{nm} > \chi_{2k-4,\beta}^2
\end{equation}
The expression for $aT_{nm}$ depends on the spectrum of $j(\hat{\alpha \mu})$ through the orbit $[U_k(\hat{\alpha \mu})]$ and the subspace spanned by $\{U_2(\hat{\alpha \mu}), \dots, U_{k-1}(\hat{\alpha \mu})\}.$ If the population antimean exists, $[U_k(\hat{\alpha \mu})]$ is a consistent estimator of $[U_k(\alpha \mu)]$, the projection on Span $\{U_2(\hat{\alpha \mu}), \dots, U_{k-1}(\hat{\alpha \mu})\}$ converges to that on Span $\{U_2(\alpha \mu), \dots, U_{k-1}(\alpha \mu)\}$ . Thus from \eqref{e:asy-chi-diff} and \eqref{e:aTnm}, $aT_{nm}$ has an asymptotic $\chi_{2k-4}^2$ distribution. Hence the test in \eqref{e:aTnm-level-beta} has asympotic level $\beta.$

\section{Application to medical imaging}

\subsection{Apert syndrome vs clinically normal children}

Our data consists of shapes for a group of eighth midface labeled anatomical landmarks from X-rays of skulls of eight year old and fourteen year-old North American children(36 boys and 26 girls), known as the {\em University School data}. {Each child's skull was imaged twice, at age 8 and next at age 14.} The data set, { collected to study anatomical changes during children growth,} represents coordinates {
 of eight craniofacial} landmarks, whose names and position on the skull are given in Bookstein ((1997)\cite{Bo:1997}), { see also
 http://life.bio.sunysb.edu/morph/data/Book-UnivSch.txt}. The data has two data set: the first one is the Apert data in Bookstein (pp. 405-406), which contains eight landmarks that describing the children who has Apert syndrome (a genetic craniosynostosis) and the second data set is clinically normal children, which contains about 40 anatomical landmarks on the skull. Out of these only 8 landmarks are registered in both groups. The two groups share only 5 registered landmarks: Posterior nasal spine,Anterior nasal spine, Sella, Sphenoethmoid registration, and Nasion. For operational definitions of these landmarks, see Bookstein (pp. 71). The shape variable of the 5 landmarks on the upper mid-face is valued in $\sum_2^5$.

\begin{figure}[H]
\centering
\begin{minipage}{.5\textwidth}
  \centering
 \includegraphics[width=80mm, height=90mm]{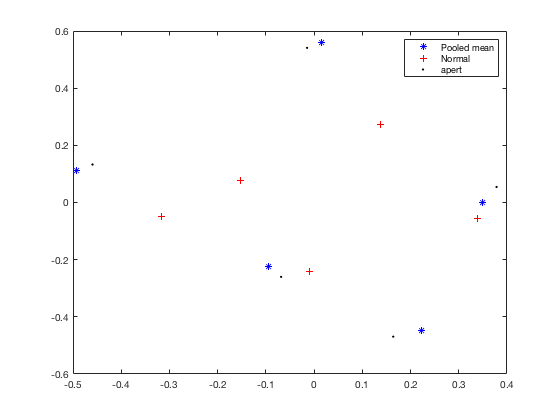}
\end{minipage}%
\begin{minipage}{.5\textwidth}
  \centering
\includegraphics[width=80mm, height=90mm]{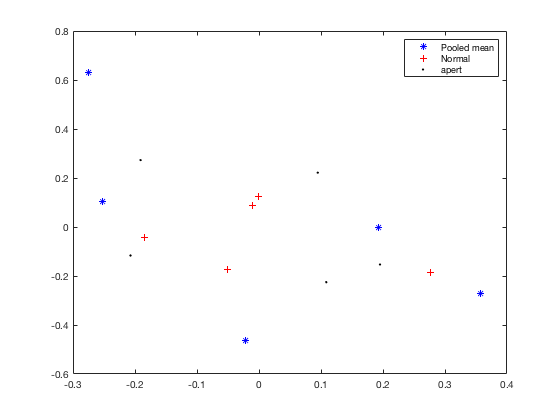}
\end{minipage}
\caption{Landmarks from preshapes of extrinsic means-left,extrinsic antimeans-right}
\end{figure}

 In our application, we take both two sample test for VW means (see Bhattacharya and Bhattacharya \cite{bhbh08}) and for VW antimeans, to see if we can distinguish between the Apert group and clinically normal group, using such Kendall shape location parameters. Figure 1 shows the plots of an icon for the sample VW means and VW antimeans for the two groups along with the pooled sample VW mean and VW antimean.

For $k=5,$ the VW mean test, the values of the two sample test statistic for testing equality of the population extrinsic mean shapes is $T_{nm}=53.1140 > \chi_{2k-4,0.05}^2=12.5916$,  along with the asymptotic $p-value=P(\chi^2_6>53.1140)=1.1129 \times 10^{-9}$ and for the VW antimean test, we get the result $aT_{nm}=144.9891 > \chi_{2k-4}^2(1-\alpha)$, along with the asymptotic $p-value=P(\chi^2_6>144.9891)=8.862936 \times 10^{-29}$.

 From this application study, we reject the null hypothesis and conclude that both the VW mean test and the VW antimean test show that one may distinguish the Apert and clinically normal group on based on their VW antimeans, or on their VW means.

\subsection{Brain scan shapes of schizophrenic vs clinically normal
children}
In this example from Bookstein (1991), 13 landmarks are recorded on a
midsagittal two-dimensional slice from a Magnetic Resonance brain scan of each of 14 schizophrenic children and 14 normal children. It is of interest to study differences in shapes of brains between the two groups which can be used to detect schizophrenia. This is an application of disease detection. The shapes of the sample k-ads lie in $\Sigma^k_2 , k = 13.$ To distinguish between the underlying distributions, we compare their VW mean and VW antimean shapes.

 For testing the equality of the VW means we use the test in Bhattacharya and Bhattacharya \cite{bhbh08}. In this application, in addition we consider the two sample test for the equality of VW antimeans developed in the previous section. Figure 2 shows the plots of the sample VW means and VW antimeans for the two groups along with the pooled sample VW mean and VW antimean for this data sets.

\begin{figure}[H]
\centering
\begin{minipage}{.5\textwidth}
  \centering
 \includegraphics[width=80mm, height=90mm]{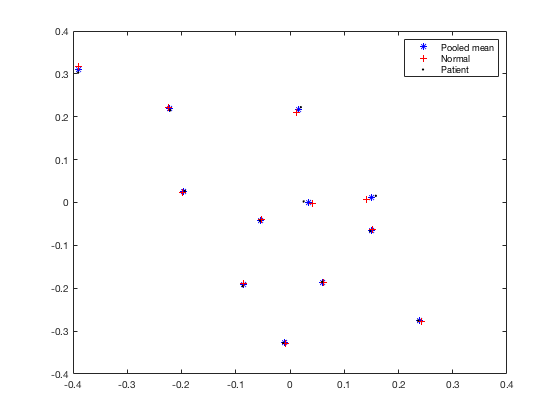}
\end{minipage}%
\begin{minipage}{.5\textwidth}
  \centering
\includegraphics[width=80mm, height=90mm]{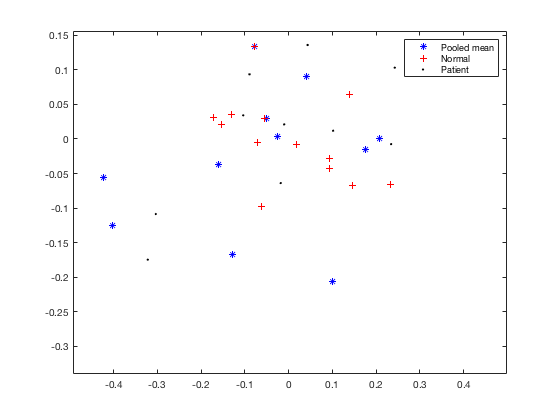}
\end{minipage}
\caption{Landmarks from preshapes of extrinsic means-left,extrinsic antimeans-right}
\end{figure}

For the VW mean test, since $k=13,$ the values of the two sample test statistic for testing equality of the population extrinsic mean shapes is $T_{nm}= 95.5476
 > \chi_{2k-4,0.05}^2=33.9244$,  along with the asymptotic $p-value=P(\chi^2_{22}>95.5476)=3.8316 \times 10^{-11}$ and for the VW antimean test, we get the result $aT_{nm}= 139.1210 > \chi_{2k-4,0.05}^2$, along with the asymptotic $p-value=P(\chi^2_{22}>144.9891)< .00001,$ which is significant at level $0.05$.

 From this application study, we therefore reject the null hypothesis and conclude that there is difference between the schizophrenic children and normal children for both the VW means or VW antimeans.
\newline

\textbf{Acknowledgments.} We would like to thank Yunfan Wang for helpful comments on an early version of our paper.

\end{document}